\documentclass{amsart}
\usepackage{amsmath,amsthm,amssymb,amscd,amsfonts}
\usepackage{enumerate}
\makeatletter
\@namedef{subjclassname@2010}{%
  \textup{2010} Mathematics Subject Classification}
\makeatother
\theoremstyle{definition}

\newcommand{\f}{\dfrac}
\newcommand{\cA}{\mathcal A}

\newtheorem{theo}{Theorem}[section]
\newtheorem{cor}[theo]{Corollary}

\newtheorem{prop}[theo]{Proposition}
\theoremstyle{definition}

\theoremstyle{definition}

\newtheorem{rem}[theo]{Remark}
\numberwithin{equation}{section}

\title{On pseudo-amenability of Beurling algebras }
\author[ K. Oustad  , A. Mahmoodi ]{ Kobra Oustad, Amin Mahmoodi}
\thanks{}
\subjclass[2010]{Primary: 22D15, 43A10; Secondary: 43A20, 46H25}
\keywords{amenability, pseudo-amenability, Beurling algebra.
 }

\address{Department of Mathematics, Central Tehran Branch,
 Islamic Azad University, Tehran, Iran.}
\address{Department of Mathematics, Central Tehran Branch,
 Islamic Azad University, Tehran, Iran.}

\email{kob.oustad.sci@iauctb.ac.ir }

\email{a\_mahmoodi@iauctb.ac.ir }

\begin{document}
\maketitle

\setcounter{section}{0}
\begin{abstract}
Amenability and pseudo-amenability of $ \ell^{1}(S,\omega) $ is
characterized, where $S$ is a left (right) zero semigroup or it is a
rectangular band semigroup. The equivalence conditions to
amenability of $\ell^{1}(S,\omega)$ are provided, where  $ S $ is a
band semigroup. For a locally compact group $G$, pseudo-amenability
of $ \ell^{1}(G,\omega) $ is also discussed.
\end{abstract}

\section{Introduction and Preliminaries }
For a Banach algebra $ \cA $ the projective tensor product $ \cA
\hat{\otimes} \cA $ is a Banach $\cA$-bimodule in a natural manner
and the multiplication map $ \pi: \cA \hat{\otimes} \cA
\longrightarrow \cA$ defined by $\pi(a \otimes b)=ab$ for $a,b \in
\cA$ is a Banach $ \cA$-bimodule homomorphism.

Amenability for Banach algebras introduced by B. E. Johnson [9]. Let
$\cA$ be a Banach algebra and $E$ be a Banach $\cA$-bimodule.  A
continuous linear operator  $D: \cA \longrightarrow E$ is a
\textit{derivation} if it satisfies $ D(ab) = D(a) \ . \ b + a \ . \
D(b) $ for all $a,b \in \cA$. Given $x \in E$, the \textit{inner}
derivation $ad_x : \cA \longrightarrow E$ is defined by $ad_x(a) = a
\ . \ x - x \ . \ a$. A Banach algebra $\cA$ is \textit{amenable} if
for every Banach $\cA$-bimodule $E$, every derivation from $\cA$
into $E^*$, the dual of $E$, is inner.

An \textit{approximate diagonal} for a Banach algebra $ \cA $ is a
net $(m_i)_i$ in  $\cA \hat{\otimes} \cA$ such that $ a \ . \ m_i -
m_i  \ . \ a \longrightarrow 0$ and $ a \pi (m_i) \longrightarrow
a$, for each $a \in \cA$. The concept of pseudo-amenability
introduced by F. Ghahramani and Y. Zhang in [5]. A Banach algebra $
\cA $ is \textit{pseudo}-\textit{amenable} if it has an approximate
diagonal. It is well-known that amenability of $ \cA $ is equivalent
to the existence of a \textit{bounded} approximate diagonal.

The notions of biprojectivity and biflatness of Banach algebras
introduced by Helemski$\breve{i}$ in [7]. A Banach algebra $ \cA $
is \textit{biprojective} if there is a bounded $ \cA $-bimodule
homomorphism $\rho: \cA \longrightarrow \cA \hat{\otimes} \cA $ such
that $ \pi o \rho=I_{\cA} $, where $I_{\cA} $ is the identity map on
$ \cA $. We say that $ \cA $ is \textit{biflat} if there is a
bounded $\cA$-bimodule homomorphism $\rho: \cA \longrightarrow (\cA
\hat{\otimes} \cA)^{**} $ such that
 $ \pi^{**} o \rho=k_{\cA} $, where $k _{\cA}: \cA \longrightarrow \cA^{**} $ is the natural embedding of $ \cA $ into its second dual.

Let $ S $ be a semigroup. A continuous function $ \omega : S
\longrightarrow (0,\infty)$ is a \textit{weight} on $S$ if  $
\omega(st) \leq \omega(s) \omega(t)$, for all $s , t \in S$. Then it
is standard that $$\ell^{1}(S,\omega)=\left\lbrace f=\sum_{s\in
S}f(s) \delta_{s}:\Vert f \Vert_{\omega}=\sum_{s\in S}\vert
f(s)\vert \omega(s)<\infty\right\rbrace  $$ is a Banach algebra with
the convolution product $ \delta_{s}\ast\delta_{t}=\delta_{st} $.
These algebras are called \textit{Beurling algebras}.

In this note, we study the earlier mentioned properties of Banach
algebras for Beurling algebras. Firstly in section 2, we
characterize amenability and pseudo-amenability of
$\ell^{1}(S,\omega)$, for some certain class of semigroups. Let $S$
be a left or right zero semigroup. We prove that pseudo-amenability
of $\ell^{1}(S,\omega)$ is equivalent to it's amenability and these
equivalent conditions imply that $S$ is singleton. We show that the
same result holds for $\ell^{1}(S,\omega)$, whenever $S$ is a
rectangular band semigroup and $ \omega$ is separable. Further, we
investigate biprojectivity of $\ell^{1}(S, \omega)$ whenever $S$ is
either left (right) zero semigroup or a rectangular band semigroup.
For a band semigroup $ S $, we show that amenability of
$\ell^{1}(S,\omega)$ is equivalent to that of $\ell^{1}(S)$ and
these are equivalent to $S$ being a finite semilattice.

Finally in section 3, we investigate pseudo-amenability of
$L^1(G,\omega)$ where $G$ is a locally compact group and $ \omega$
is a weight on $G$. We prove that pseudo-amenability of
$L^1(G,\omega)$ implies amenability of $G$, and under a certain
condition it implies diagonally boundedness of $ \omega$. Next, if
$L^1(G,\omega)$ is pseudo-amenable we may obtain a character $
\varphi$ on $G$ for which $\varphi \leq \omega$.

\section{Amenability and pseudo-amenability of $ \ell^{1}(S,\omega) $}
 A semigroup $ S $ is a \textit{left zero
semigroup } if $ st=s$, and it is a \textit{right zero semigroup }
if $ st=t$ for each $ s,t\in S $. Then for $
f,g\in\ell^{1}(S,\omega) $, it is obvious that $f * g =
\varphi_{S}(f)g$ if $S$ is a right zero semigroup, and $f * g =
\varphi_{S}(g)f$ if $S$ is a left zero semigroup, where $\varphi_{S}
$ is the \textit{augmentation character} on $ \ell^{1}(S,\omega) $.

We extend somewhat the obtained results for $\ell^{1}(S)$ in [2,3]
to the weighted case $\ell^{1}(S,\omega)$.
\begin{prop} Suppose that $ S $ is a right (left) zero semigroup and $ \omega $ be a weight
on $ S $. Then $ \ell^{1}(S,\omega) $ is biprojective.
\end{prop}
\textbf{Proof.} We only give the proof in the case $S$ is a right
zero semigroup. Define $\rho: \ell^{1}(S,\omega) \longrightarrow
\ell^{1}(S,\omega) \otimes \ell^{1}(S,\omega)$ by $\rho(f)=
\delta_{t_{\circ}} \otimes f,$, where $t_0$ is an arbitrary element
$S$. Then for each $ f,g\in\ell^{1}(S,\omega) $ we have $$\rho(f
\ast g)=\delta_{t_{\circ}} \otimes (f \ast g)= \varphi_{S}(f)
(\delta_{t_{\circ}} \otimes g) = ( f * \delta_{t_0}) \otimes g = f \
. \ (\delta_{t_0} \otimes g) = f \ . \ \rho (g)$$ and similarly
$\rho (f \ast g)= \rho(f) \ . \ g.$ Further, $\pi \rho$ is the
identity map on  $ \ell^{1}(S,\omega) $, as required. \qed

\begin{rem}
It is known that every biprojective Banach algebra is  biflat. Hence
Proposition 2.1 shows that for every right or left zero semigroup
$S$,  $ \ell^{1}(S,\omega) $ is biflat.
\end{rem}
Given two semigroups $S_1$ and $S_2$, we say that a weight $\omega$
on $ S:= S_1 \times S_2$ is \textit{separable} if there exist two
weights $ \omega_{1} $ and $ \omega_{2} $ on  $S_1$ and $S_2$,
respectively such that $ \omega=\omega_{1}\otimes \omega_{2} $. It
is easy to verify that $  \ell^{1}(S,\omega) \cong
\ell^{1}(S_1,\omega_1) \hat{\otimes} \ell^{1}(S_2,\omega_2)$.

Let $ S $ be a semigroup and let $ E(S)=\lbrace p\in
S:p^{2}=p\rbrace $. We say that $ S $ is a \textit{band semigroup}
if $ S=E(S) $. A band semigroup $S$ satisfying $sts=s$, for each $s,
t \in S$ is called a \textit{rectangular band semigroup}. For a
rectangular band semigroup $S$, it is known that $ S \simeq L \times
R $, where $ L $ and $ R $ are left and right zero semigroups,
respectively [8, Theorem 1.1.3].

\begin{prop} Let $ S $ be a rectangular band semigroup and $\omega $ be a
separable weight on $ S $. Then $ \ell^{1}(S,\omega) $ is
biprojective, and so it is biflat.
\end{prop}
\textbf{Proof.} In view of earlier argument, it follows From
Proposition 2.1, and then from [10, Proposition 2.4].

\begin{theo}
Let $ S $ be a rectangular band semigroup and $\omega $ be a weight
on $ S $. Then $ \ell^{1}(S,\omega) $ is amenable if and only if $ S
$ singleton.
\end{theo}
\textbf{Proof.} From [11, Theorem 3.6], $ \ell^{1}(S) $ is amenable.
Then it is immediate by [2, Theorem 3.3]. \qed

For a semigroup $S$, we denote by $S^{op}$ the semigroup whose
underlying space is $S$ but whose multiplication is the
multiplication in $S$ reversed.

\begin{prop} Let $ S $ be a right (left) zero semigroup and $ \omega $ be a weight
on $ S $. Then $ \ell^{1}(S,\omega) $ is amenable if and only if $S$
is singleton.
\end{prop}
{\bf Proof.} Suppose that $ S $ is a left zero semigroup, and that $
\ell^{1}(S,\omega) $ is amenable. Then $S^{op}$ is a right zero
semigroup. It is readily seen that $S \times S^{op}$ is a
rectangular band semigroup, and $ \ell^{1}(S^{op},\omega) $ is
amenable. Hence $ \ell^{1}(S,\omega) \hat{\otimes}
\ell^{1}(S^{op},\omega) \cong \ell^{1}(S \times S^{op} , \omega
\otimes \omega)$ is amenable. Now, by Theorem 2.4, $S$ is singleton.
\qed

Let $ \cA $ be Banach algebra, $ \mathcal{I} $ be a
\textit{semilattice} (i.e., $ \mathcal{I} $ is a commutative band
semigroup) and $\{ \cA _{\alpha}:\alpha \in \mathcal{I} \}$ be a
collection of closed subalgebras of $\cA $. Then $\cA $ is $
\ell^{1}$-\textit{graded} of $ \cA _{\alpha}$'s over the semilattice
$ \mathcal{I} $, denoted by $ \cA = \bigoplus_{\alpha \in
\mathcal{I}}^{\ell^{1}} \cA_{\alpha} $, if it is $
\ell^{1}$-directsum of $ \cA _{\alpha}$'s as Banach space such that
$\cA _{\alpha} \cA _{\beta} \subseteq \cA _{\alpha \beta}$, for each
$ \alpha , \beta \in \mathcal{I}$.

Suppose that $ S^{1} $ is the unitization of a semigroup $ S $. An
equivalence relation $ \tau $ on $ S $ is defined by $s \tau t
\Longleftrightarrow S^{1} s S^{1}=S^{1} t S^{1}$, for all $s, t \in
S$. If  $ S $ is a band semigroup, then by [8, Theorem 4.4.1], $ S=
\bigcup_{\alpha\in \mathcal{I}} S_{\alpha} $ is a semilattice of
rectangular band semigroups, where $ \mathcal{I} = \frac{S}{\tau} $
and for each $ \alpha=[s] \in \mathcal{I}$, $ S_{\alpha}=[s] $.

\begin{theo} Let $ S $ be a band semigroup and $ \omega $ be a weight on $ S $.
Then the following are equivalent:

$(i)$  $ \ell^{1}(S,\omega) $ is amenable.

$(ii)$ $ S $ is finite and each $\tau-$class is singleton.

$(iii)$  $ \ell^{1}(S) $ is amenable.

$(iv)$ $ S $ is a finite semilattice.
\end{theo}
\textbf{Proof:} The implications $(ii)$ to $(iv)$ are equivalent [2,
Theorem 3.5]. We establish $(i) \longrightarrow  (ii) $ and  $(iv)
\longrightarrow  (i) $.

$(i) \longrightarrow  (ii) $ If $ \ell^{1}(S,\omega) $ is amenable,
then $ E(S)=S $ is finite and so $ \mathcal{I}= \frac{S}{\tau} $ is
a finite semilattice. Hence $\ell^{1}(S,\omega) \cong
\bigoplus_{\alpha \in \mathcal{I}}^{\ell^{1}}
\ell^{1}(S_{\alpha},\omega_{\alpha})$, where $ \omega_{\alpha} =
\omega \vert_{S_{\alpha}} $. Then by [6, Proposition 3.1], each $
\ell^{1}(S_{\alpha},\omega_{\alpha}) $ is amenable. Now by Theorem
2.4, $ S_{\alpha} $ is singleton for each $ \alpha \in \mathcal{I}
$, as required.

$(iv) \longrightarrow  (i) $ In this case  $\ell^{1}(S,\omega) \cong
 \ell^{1}(S) $, and $\ell^{1}(S) $ is amenable. \qed

\begin{theo} Let $ S $ be a rectangular band semigroup, and let $ \omega $ be a separable weight on $ S $.
Then $ \ell^{1}(S,\omega) $ is pseudo-amenable if and only if $S$ is
singleton.
\end{theo}
\textbf{Proof.}  There is a left zero semigroup $L$ and a right zero
semigroup $R$, and there are weights $ \omega_L$ and $ \omega_R$ on
$L$ and $R$, respectively such that $ S \cong L \times R$ and $
\omega = \omega_L \otimes \omega_R$. We have $ \ell^{1}(S,\omega)
\cong \ell^{1}(L, \omega_L) \hat{\otimes} \ell^{1}( R, \omega_R)$.
Hence the map  $ \theta: \ell^{1}(S,\omega) \longrightarrow
\ell^{1}(L, \omega_L)$ defined by $ \theta(f \otimes
g)=\varphi_{R}(g)f$ for $ f \in \ell^{1}(L, \omega_L)$ and $g \in
\ell^{1}( R, \omega_R) $, is an epimorphism of Banach algebras,
whereas $ \varphi_{R} $ is the \textit{augmentation} character on $
\ell^{1}( R, \omega_R) $. Whence $ \ell^{1}(L, \omega_L) $ has left
and right approximate identity. Therefore $ L $ is singleton,
because it is left zero semigroup. Similarly $ R $ is singleton, so
is $ S $. \qed

\begin{cor} Let $ S $ be a right (left) zero
semigroup and $ \omega $ be a weight on $ S $. Then the following
are equivalent:

$(i)$ $ \ell^{1}(S,\omega) $ is pseudo-amenable.

$(ii)$  $S$ is singleton.

$(iii)$ $ \ell^{1}(S,\omega) $ is amenable.
\end{cor}
{\bf Proof.} The implication $(ii) \longleftrightarrow (iii)$ is
Proposition 2.5. For $(i)\longrightarrow (ii)$, we apply Theorem 2.7
for the rectangular band semigroup $S \times S^{op} $ with $
\omega_L = \omega_R = \omega$. \qed

The following is a combination of Theorems 2.4 and 2.7. Notice that
in Theorem 2.4,  we need not $ \omega $ to be separable.
\begin{cor} Let $ S $ be a rectangular band semigroup, and let $ \omega $ be a separable weight on $ S
$. Then the following are equivalent:

$(i)$ $ \ell^{1}(S,\omega) $ is pseudo-amenable.

$(ii)$  $S$ is singleton.

$(iii)$ $ \ell^{1}(S,\omega) $ is amenable.
\end{cor}
For the left cancellative semigroups we have the following.
\begin{theo}
Suppose that $ S $ is a left cancellative semigroup and $\omega $ is
a weight on $S$. If $ \ell^{1}(S,\omega) $ is pseudo-amenable, then
$ S $ is a group.
\end{theo}
\textbf{proof:} This is a more or less verbatim of the proof of [3,
Theorem 3.6 $(i) \longrightarrow (ii)$]. \qed
\section{ Pseudo-amenability of $L^1(G,\omega)$}

Throughout $G$ is a locally compact group and $\omega$ is a weight
on $G$. The weight $\omega$ is \textit{diagonally bounded} if $
\sup_{g \in G} \omega(g) \omega(g^{-1}) < \infty$. It seems to be a
\textit{right} conjecture that $L^1(G,\omega)$ will fail to be
pseudo-amenable whenever $\omega$ is not diagonally bounded.
Although we are not able to prove (or disprove) the conjecture, the
following is a weaker result.

The proofs in this section owe much to those of [4, Section 8].
\begin{theo} Suppose that $L^1(G,\omega)$ is pseudo-amenable for which there is an approximate diagonal $(m_i)_i$ such that
 $
m_i - \delta_g \ . \ m_i \ . \ \delta_{g^{-1}} \longrightarrow 0$
uniformly on $G$. Then $\omega$ is diagonally bounded.
\end{theo}
{\bf Proof.} We follow the standard argument in [4, Proposition
8.7]. Choose $f \in L^1(G,\omega)$ such that $K:= supp f$ is compact
and $\int f \neq 0$. Putting $ F:= f \ . \ \chi_K \in L^\infty (G,
\omega^{-1})$, we see that $\pi^*(F) \in L^\infty (G \times G ,
\omega^{-1} \times \omega^{-1})$ with $$ \pi^*(F) (x,y) = F(xy) =
\int \chi_K (xyt) f(t) dt \ .$$ Let $(m_i)_i \subseteq L^1(G \times
G ,\omega \times \omega)$ be an approximate diagonal for
$L^1(G,\omega)$ such that $\delta_g \ . \ m_i \ . \ \delta_{g^{-1}}
- m_i \longrightarrow 0$ uniformly on $G$, and $ \pi(m_i)  f - f
\longrightarrow 0$. Then for each $i$ $$\langle \pi^*(F) , m_i
\rangle = \langle F , \pi(m_i) \rangle = \langle \chi_K , \pi(m_i) f
\rangle \longrightarrow \langle \chi_K ,  f \rangle = \int f \ .$$
Consequently
$$ \lim_i \langle \pi^*(F) , m_i \rangle \neq 0 \ . \ \ \ \ (1)$$ We
define $E:= K K^{-1}$, and $ A := \{ (x,y) \in G \times G \ : \ xy
\in E \}.$ For $r > 0$, we define $A_r := \{ (x,y) \in A \ : \
\omega(x) \omega(y) < r \},$ and $B_r := \{ (x,y) \in A \ : \
\omega(x) \omega(y) \geq r \}$. Obviously, $ \pi^*(F) \chi_{A_r}$
and $ \pi^*(F) \chi_{B_r}$ both are in $L^\infty (G \times G ,
\omega^{-1} \times \omega^{-1})$, and $ \pi^*(F) = \pi^*(F) \chi_A =
\pi^*(F) \chi_{A_r} +  \pi^*(F) \chi_{B_r}$. For every $i$, it is
easy to see that $$ | \langle \pi^*(F) \chi_{B_r} , m_i    \rangle |
\leq ||m_i|| \  ||F|| \  r^{-1}  \ c_1
$$ where $c_1 := \sup_{t \in E} \omega(t)$. Hence $$ \lim_{r
\longrightarrow \infty}  \langle  \pi^*(F) \chi_{B_r} , m_i \rangle
= 0 \ . \ \ \  \ (2) $$ Next, for every $g \in G$, $r > 0$, and $i$,
we obtain $$ | \langle  \pi^*(F) \chi_{A_r} ,  \delta_g \ . \ m_i \
. \  \delta_{g^{-1}}    \rangle |  \leq || m_i|| \ ||F|| \  r \  c_1
 \  c_2^2 \ \f{1}{\omega(g) \omega(g^{-1})}$$ where $c_2:=\sup_{t \in
E^{-1}} \omega(t)$. Therefore
\begin{align*}
| \langle  \pi^*(F) \chi_{A_r} ,  m_i    \rangle | & \leq | \langle
\pi^*(F) \chi_{A_r} ,  m_i - \delta_g \ . \ m_i \ . \
\delta_{g^{-1}}    \rangle | + | \langle \pi^*(F) \chi_{A_r} ,
\delta_g \ . \ m_i \ . \ \delta_{g^{-1}}    \rangle | \\& \leq ||
\pi^*(F) || \ \sup_{g \in G} || m_i - \delta_g \ . \ m_i \ . \
\delta_{g^{-1}} ||    + || m_i|| \ ||F|| \  r \ c_1 \ c_2^2 \
\f{1}{\omega(g) \omega(g^{-1})} \ . \ \ \ \ (3)
\end{align*}

Towards a contradiction, we assume that $\omega$ is not diagonally
bounded. Then there is a sequence $(g_n)_n$ in $G$ such that $\lim_n
\omega(g_n) \omega(g_n^{-1}) = \infty \ .$  Whence, it follows from
$(3)$ that for each $i$ and $r > 0$ $$| \langle  \pi^*(F) \chi_{A_r}
, m_i \rangle | \leq
 ||
\pi^*(F) || \  \sup_{g \in G} || m_i - \delta_g \ . \ m_i \ . \
\delta_{g^{-1}} ||  \ . \ \ \ \ (4) $$ Hence
$$ | \langle \pi^*(F)  ,  m_i \rangle | \leq || \pi^*(F)
|| \ \sup_{g \in G} || m_i - \delta_g \ . \ m_i \ . \
\delta_{g^{-1}} || + | \langle \pi^*(F) \chi_{B_r} , m_i \rangle | \
.$$ Putting $(2)$ and $(4)$ together, we may see that $$ \lim_i
\langle \pi^*(F) , m_i \rangle  = 0$$ contradicting $(1)$. \qed

\begin{theo} Suppose that $L^1(G,\omega)$ is pseudo-amenable, and that $\omega$ is bounded away from 0. Then
$G$ is amenable.
\end{theo}
{\bf Proof.} Since  $L^1(G,\omega)$ is unital, pseudo-amenability
and approximate amenability are the same [5, Proposition 3.2]. Now,
it is immediate by [4, Proposition 8.1]. \qed

We conclude by the following which is an analogue of [4, Proposition
8.9].
\begin{prop} Let $L^1(G,\omega)$ be pseudo-amenable. Then there is a
continuous positive character $\varphi$ on $G$ such that $\varphi
\leq \omega$.
\end{prop}
{\bf Proof.} Suppose that $(m_i)_i \subseteq L^1(G \times G ,\omega
\times \omega)$ be an approximate diagonal for $L^1(G,\omega)$. For
each $i$ and $f \in L^\infty (G \times G, \omega^{-1} \times
\omega^{-1})^+$ we define $$ \widetilde{m_i}(f):= \sup \{ Re \langle
m_i , \psi \rangle \ : \ 0 \leq | \psi | \leq f , \ \psi \in
L^\infty (G \times G, \omega^{-1} \times \omega^{-1}) \}.$$ Then
$\widetilde{m_i} \neq 0$ on $L^\infty (G \times G, \omega^{-1}
\times \omega^{-1})^+$ and we may extend $\widetilde{m_i}$ to a
bounded linear functional on $L^\infty (G \times G, \omega^{-1}
\times \omega^{-1})$ in the obvious manner. It is readily seen that
$\widetilde{m_i} \neq 0$, $\langle \widetilde{m_i} , f \rangle \geq
0$, and $\delta_{g^{-1}} \ . \ \widetilde{m_i} \ . \ \delta_g -
\widetilde{m_i} \longrightarrow 0$, for every $f \in L^\infty (G
\times G, \omega^{-1} \times \omega^{-1})^+$ and $g \in G$.

Putting $\widetilde{\omega} (x) := \sup_{g \in G} \omega(g^{-1} x
g)$, $x \in G$. Then $\widetilde{\omega}  \in L^\infty (G,
\omega^{-1})$, $\widetilde{\omega} (xy) = \widetilde{\omega} (yx)$,
$ \pi^*(\widetilde{\omega}) \in L^\infty (G \times G, \omega^{-1}
\times \omega^{-1})$, and $ \delta_g \ . \ \pi^*(\widetilde{\omega})
\ . \ \delta_{g^{-1}} = \pi^*(\widetilde{\omega})$.

Take $f \in C_c(G)^+$ with $\int f = 1$, and then $h:= f \ . \
\chi_K$, where $K:= supp f$. One may see that $h$ is continuous, and
there is $c > 0$ such that $ \pi^*(\widetilde{\omega}) \geq c
\pi^*(h)$. Hence
\begin{align*}
\lim_i \langle \widetilde{m_i} , \pi^*(\widetilde{\omega}) \rangle &
\geq c \lim_i \langle \widetilde{m_i} , \pi^*(h) \rangle \geq c
\lim_i Re \langle m_i , \pi^*(h) \rangle = c \lim_i Re \langle
\pi(m_i) , h \rangle \\& = c \lim_i Re \langle \pi(m_i) \ . \ f ,
\chi_K \rangle = c Re \langle f , \chi_K \rangle = c > 0 \ .
\end{align*}
Therefore there is $i_0$ for which $ \langle \widetilde{m_{i_0}} ,
\pi^*(\widetilde{\omega}) \rangle > 0$. Set $F:= \langle
\widetilde{m_{i_0}} , \pi^*(\widetilde{\omega}) \rangle^{-1}
\pi^*(\widetilde{\omega})$, and for $g \in G$ we put $$  A_g (x,y)
:= \f{1}{2} ( \log \f{\omega(gx) \omega(gy^{-1}) }{\omega(x)
\omega(y^{-1}) } ) F(x,y) \ , \ \ \ (x,y \in G) \ . $$ Finally, for
each $g \in G$, we define $ \varphi(g) := \exp \langle
\widetilde{m_{i_0}} , A_g \rangle. $ A similar argument used in [4,
Proposition 8.9], shows that $\varphi$ is the desired character on
$G$. \qed

\end{document}